\documentclass[a4paper,11pt,twoside,reqno]{amsart}

\usepackage[utf8]{inputenc}
\usepackage[plainpages=false,pdfpagelabels=true]{hyperref}
\usepackage{amssymb,amsthm,slashed,dsfont}
\usepackage[margin=1in]{geometry}
\usepackage[mathscr]{eucal}

\newtheorem{Satz}{Theorem}[section]
\newtheorem{Prop}[Satz]{Proposition}
\newtheorem{Lem}[Satz]{Lemma}

\theoremstyle{definition}

\newtheorem{Bem}[Satz]{Remark}

\newcommand{\tr}{\operatorname{tr}}

\newcommand{\dv}{\text{ }dV}

\allowdisplaybreaks[1]

\renewcommand{\epsilon}{\varepsilon}

\newcommand{\R}{\ensuremath{\mathbb{R}}}

\numberwithin{equation}{section}

\title{The stress-energy tensor for polyharmonic maps}
\author{Volker Branding}
\date{\today}
\address{University of Vienna, Faculty of Mathematics\\
Oskar-Morgenstern-Platz 1, 1090 Vienna, Austria\\}
\email{volker.branding@univie.ac.at}
\begin{document}

\begin{abstract}
We derive the stress-energy tensor for polyharmonic maps between Riemannian manifolds.
Moreover, we employ the stress-energy tensor to characterize polyharmonic maps where
we pay special attention to triharmonic maps.
\end{abstract} 

\maketitle

\section{Introduction and results}
Harmonic maps are one of the most studied geometric variational problems in mathematics.
Besides their importance in analysis and geometry they also appear frequently in the physics literature
for example as nonlinear \(\sigma\)-model in quantum field theory.

In order to define harmonic maps we consider two Riemannian manifolds \((M,g)\) and \((N,h)\).
Given a smooth map \(\phi\colon M\to N\) we can define its energy by
\begin{align*}
E(\phi)=\int_M|d\phi|^2\dv.
\end{align*}
Its critical points are characterized by the vanishing of the \emph{tension field} which is given by
\begin{align}
\label{harmonic-map-equation}
0=\tau(\phi):=\tr_g\nabla d\phi=0,\qquad \tau(\phi)\in\Gamma(\phi^\ast TN).
\end{align}
Solutions of \eqref{harmonic-map-equation} are called \emph{harmonic maps}. 
The harmonic map equation constitutes a semilinear second order partial differential equation.
There is a vast number of results on existence and qualitative behavior of harmonic maps, we refer to \cite{Branding2019jgea} and references therein for an overview.

A generalization of harmonic maps that receives growing attention is 
given by the theory of \emph{biharmonic maps}. These are critical points of the \emph{bienergy}
\begin{align*}
E_2(\phi)=\int_M|\tau(\phi)|^2\dv
\end{align*}
and are characterized by the vanishing of the bitension field 
\begin{align}
\label{biharmonic-map-equation}
0=\tau_2(\phi):=-\Delta\tau(\phi)-\tr_g R^N(d\phi(\cdot),\tau(\phi))d\phi(\cdot).
\end{align}
Here, \(\Delta\) denotes the rough Laplacian on the vector bundle \(\phi^\ast TN\) and \(R^N\)
represents the curvature tensor on the target \(N\).
The biharmonic map equation \eqref{biharmonic-map-equation} is a semilinear elliptic partial differential equation of fourth order
which makes its analysis more involved compared to the case of harmonic maps.
It is obvious that every harmonic map is also biharmonic, on the other hand 
a biharmonic map can be non-harmonic in which case it is called \emph{proper biharmonic}.
Many analytic and geometric conditions are known that force a biharmonic map to be harmonic,
see for example \cite{MR3834926,branding2018nonexistence} and references therein for a recent overview.

One systematic generalization of both harmonic and biharmonic maps is given by the so-called
\emph{polyharmonic maps}.
To define the latter one can study the critical points 
of the following energy functionals, where one has to distinguish between polyharmonic 
maps of even and odd order.
In the even case we set
\begin{align}
\label{poly-energy-even}
E_{2s}(\phi)=\int_M|\Delta^{s-1}\tau(\phi)|^2\dv,
\end{align}
whereas in the odd case we have
\begin{align}
\label{poly-energy-odd}
E_{2s+1}(\phi)=\int_M|\nabla\Delta^{s-1}\tau(\phi)|^2\dv.
\end{align}
Here, \(s=1,2,\ldots\) is a natural number.

The first variations of \eqref{poly-energy-even} and \eqref{poly-energy-odd} were calculated in \cite{MR3007953},
its critical points are referred to as polyharmonic maps (of order \(k\)) or \(k\)-harmonic maps.
We have to distinguish two cases:
\begin{enumerate}
 \item If \(k=2s,s=1,2,\ldots\), the critical points of \eqref{poly-energy-even} are given by
\begin{align}
\label{tension-2s}
0=\tau_{2s}(\phi):=&\Delta^{2s-1}\tau(\phi)-R^N(\Delta^{2s-2}\tau(\phi),d\phi(e_j))d\phi(e_j) \\
 \nonumber&+\sum_{l=1}^{s-1}\big(R^N(\Delta^{s-l-1}\tau(\phi),\nabla_j\Delta^{s+l-2}\tau(\phi))d\phi(e_j) \\
\nonumber&\hspace{0.5cm}-R^N(\nabla_j\Delta^{s-l-1}\tau(\phi),\Delta^{s+l-2}\tau(\phi))d\phi(e_j)\big).
\end{align}
\item If \(k=2s+1,s=0,1,\ldots\), the critical points of \eqref{poly-energy-odd} are given by
\begin{align}
\label{tension-2s+1}
0=\tau_{2s+1}(\phi):=&\Delta^{2s}\tau(\phi)-R^N(\Delta^{2s-1}\tau(\phi),d\phi(e_j))d\phi(e_j) \\
\nonumber&-\sum_{l=1}^{s-1}\big(R^N(\nabla_j\Delta^{s+l-1}\tau(\phi),\Delta^{s-l-1}\tau(\phi))d\phi(e_j) \\
\nonumber&\hspace{0.5cm}-R^N(\Delta^{s+l-1}\tau(\phi),\nabla_j\Delta^{s-l-1}\tau(\phi))d\phi(e_j)\big) \\
\nonumber&\hspace{0.5cm}-R^N(\nabla_j\Delta^{s-1}\tau(\phi),\Delta^{s-1}\tau(\phi))d\phi(e_j).
\end{align}
\end{enumerate}
Here, we have set \(\Delta^{-1}=0\). In addition, \(\{e_j\},j=1,\ldots,\dim M=m\) denotes an orthonormal basis of \(TM\)
and we are applying the Einstein summation convention. Moreover, we use the short notation \(\nabla_i:=\nabla_{e_i}\).

It can be directly seen that every harmonic map is a solution of the polyharmonic map equations \eqref{tension-2s} and \eqref{tension-2s+1},
see \cite{branding2019structure} for a recent classification result.

In this article we focus on the stress-energy tensor associated to both \eqref{poly-energy-even} and \eqref{poly-energy-odd}.
The stress-energy tensor of a given energy functional can be derived by varying the energy functional
with respect to the domain metric. The stress-energy tensor calculated in this fashion is divergence-free
assuming that we are considering a critical point of the energy functional.

This type of conservation law is a consequence of Noether's theorem as the energy functionals \eqref{poly-energy-even}
and \eqref{poly-energy-odd} are invariant under diffeomorphisms on the domain which we will prove later.

We want to mention that in the physics literature one usually uses the phrase \emph{energy-momentum tensor}
instead of stress-energy tensor.

For harmonic maps the stress-energy tensor was calculated in \cite{MR655417},
for biharmonic maps it was stated in \cite{MR891928} and later systematically derived in \cite{MR2395125}.
In this article we will extend these calculations to the case of polyharmonic maps.

Throughout this article we will use the following sign conventions.
For the Riemannian curvature tensor we use \(R(X,Y)Z=[\nabla_X,\nabla_Y]Z-\nabla_{[X,Y]}Z\) and
for the (rough) Laplacian on \(\phi^\ast TN\) we use the geometer's sign convention \(\Delta:=-\tr_g(\nabla\nabla-\nabla_\nabla)\).

We will use the same symbol $\langle\cdot,\cdot\rangle$ to denote the Riemannian metrics on various vector bundles,
and the same symbol $\nabla$ for the corresponding Riemannian connections. Whenever necessary we will be more precise
and explicitly denote the connection in use.

Whenever we will make use of indices, we will use
Latin indices \(i,j,k\) for indices on the domain ranging from \(1\) to \(m\)
and Greek indices \(\alpha,\beta,\gamma\) for indices on the target
which take values between \(1\) and \(n\). In addition, local coordinates on the domain
will be denoted by \(x^i,i=1,\ldots,\dim M=m\) and for local coordinates on the the target
we will use \(y^\alpha,\alpha=1,\ldots,\dim N=n\).

We will employ the Einstein summation convention, that is we will sum over repeated indices.

This article is organized as follows: In Section 2 we derive the 
stress-energy tensor for polyharmonic maps and discuss the invariance of the energy functionals
\eqref{poly-energy-even} and \eqref{poly-energy-odd} under diffeomorphisms on the domain.
In Section 3 we study maps with vanishing stress-energy tensor. Finally, in Section 4 we use
the stress-energy tensor to obtain a non-existence result for
triharmonic maps with finite energy.

\section{Deriving the stress-energy tensor for polyharmonic maps}
We derive the stress-energy tensor associated to the energy functionals \eqref{poly-energy-even}
and \eqref{poly-energy-odd} by varying the functionals with respect to the metric on the domain. 
To this end we set
\begin{align}
\label{variation-metric-domain}
\frac{d}{dt}\big|_{t=0}g_{ij}=\omega_{ij},
\end{align}
where \(\omega_{ij}\) is a smooth symmetric \(2\)-tensor on \(M\).

We have to distinguish between the cases of polyharmonic maps of even and odd order.
\subsection{The even case}
First, suppose that we have a polyharmonic map of even order.
\begin{Lem}
\label{even-lem-em-a}
Let \(\phi\colon M\to N\) be a smooth map and consider a variation of the metric on \(M\) as defined in \eqref{variation-metric-domain}.
Then the following identity holds
\begin{align}
\frac{d}{dt}\big|_{t=0}E_{2s}(\phi,g_t)=&\frac{1}{2}\int_M|\Delta^{s-1}\tau(\phi)|^2\langle g,\omega\rangle\dv+2\int_M\langle\frac{d}{dt}\big|_{t=0}\Delta^{s-1}\tau(\phi),\Delta^{s-1}\tau(\phi)\rangle\dv.
\end{align}
\end{Lem}

\begin{proof}
The variation of the volume element\(\dv\) is given by
\begin{align*}
\frac{d}{dt}|_{t=0}\dv_{g_t}=\frac{1}{2}\langle g,\omega\rangle\dv_g.
\end{align*}
The claim then follows by a direct calculation.
\end{proof}

In the following we require the variation of the tension field of the map \(\phi\) 
with respect to the domain metric \(g\), which was already calculated in \cite{MR2395125}.
We will perform this calculation in local coordinates.
To this end let \((U,x^i)\) be a local coordinate chart on \(M\) and \((V,y^\alpha)\) a local coordinate chart on \(N\)
such that \(\phi(U)\subset V\).

\begin{Lem}
\label{lem-variation-tension}
Let \(\phi\colon M\to N\) be a smooth map and consider a variation of the metric on \(M\) as defined in \eqref{variation-metric-domain}.
The variation of the tension field with respect to the metric on the domain is given by
\begin{align}
\label{variation-tension-field}
\frac{d}{dt}\big|_{t=0}\tau^\alpha(\phi)=-\omega^{ij}(\nabla d\phi)^\alpha_{ij}-(\nabla_i\omega^{ki})\frac{\partial\phi^\alpha}{\partial x^k}
+\frac{1}{2}(\nabla^k\tr\omega)\frac{\partial\phi^\alpha}{\partial x^k},
\end{align}
where \(\alpha=1,\ldots,n\).

The variation of the connection Laplacian on \(\phi^\ast TN\) with respect to the metric on the domain is given by
\begin{align}
\label{variation-connection-laplacian}
\frac{d}{dt}\big|_{t=0}\Delta=\omega^{ij}\nabla_i\nabla_j+(\nabla^i\omega^j_{~i})\nabla_j-\frac{1}{2}(\nabla^k\tr\omega)\nabla_k.
\end{align}
\end{Lem}

\begin{proof}
In terms of local coordinates we may express the tension field as
\begin{align*}
\tau^\alpha(\phi)=g^{ij}\big(\frac{\partial^2\phi^\alpha}{\partial x^i\partial x^j}
-\Gamma^k_{ij}\frac{\partial\phi^\alpha}{\partial x^k}
+\Gamma^\alpha_{\beta\gamma}\frac{\partial\phi^\beta}{\partial x^i}\frac{\partial\phi^\gamma}{\partial x^j}
\big),
\end{align*}
where \(\Gamma^k_{ij}\) are the Christoffel symbols on \(M\) and \(\Gamma^\alpha_{\beta\gamma}\) 
are the Christoffel symbols on \(N\).

Now, consider a variation of the metric on \(M\) as defined in \eqref{variation-metric-domain}.
Then the following formula holds
\begin{align*}
\frac{d}{dt}\big|_{t=0}\Gamma^k_{ij}&=\frac{1}{2}g^{kr}(\nabla_j\omega_{ri}+\nabla_i\omega_{rj}-\nabla_r\omega_{ij}),
\end{align*}
see for example \cite[Equation (2.5)]{MR2395125}.

Note that \eqref{variation-metric-domain} directly implies that
\begin{align*}
\frac{d}{dt}\big|_{t=0}g^{ij}=-\omega^{ij}.
\end{align*}

Moreover, we find using the variation of the connection on \(M\) that
\begin{align*}
g^{ij}\big(\frac{d}{dt}\big|_{t=0}\Gamma^k_{ij}\big)=\nabla_i\omega^{ki}-\frac{1}{2}\nabla^k \tr\omega,
\end{align*}
which already yields the first claim.

The connection Laplacian on the vector bundle \(\phi^\ast TN\) has the following local expression
\begin{align*}
\Delta=-g^{ij}(\nabla_i\nabla_j-\Gamma^k_{ij}\nabla_k)
\end{align*}
such that its variation with respect to the metric on \(M\) is given by
\begin{align*}
\frac{d}{dt}\big|_{t=0}\Delta=-\big(\frac{d}{dt}g^{ij}\big|_{t=0}\big)(\nabla_i\nabla_j-\Gamma^k_{ij}\nabla_k)
+g^{ij}\big(\frac{d}{dt}\Gamma^k_{ij}\big|_{t=0}\big)\nabla_k,
\end{align*}
which completes the proof.
\end{proof}

\begin{Lem}
Let \(\phi\colon M\to N\) be a smooth map and consider a variation of the metric on \(M\) as defined in \eqref{variation-metric-domain}.
Then the following identity holds
\begin{align}
\label{even-lem-em-c}
\int_M\langle\frac{d}{dt}\big|_{t=0}\Delta^{s-1}\tau(\phi),\Delta^{s-1}\tau(\phi)\rangle\dv
=&
\sum_{l=1}^{s-1}\int_M\langle\big(\frac{d}{dt}\big|_{t=0}\Delta\big)\Delta^{s-l-1}\tau(\phi),\Delta^{s+l-2}\tau(\phi)\rangle\dv \\
\nonumber&+\int_M\langle\frac{d}{dt}\big|_{t=0}\tau(\phi),\Delta^{2s-2}\tau(\phi)\rangle\dv.
\end{align}
\end{Lem}

\begin{proof}
First of all, we note that
\begin{align*}
\frac{d}{dt}\big|_{t=0}\Delta^{s-1}\tau(\phi)=&\sum_{l=1}^{s-1}\Delta^{l-1}\big(\frac{d}{dt}\big|_{t=0}\Delta\big)\Delta^{s-l-1}\tau(\phi)
+\Delta^{s-1}\frac{d}{dt}\big|_{t=0}\tau(\phi). 
\end{align*}
The claim then follows by using the self-adjointness of the Laplacian in \(L^2\).
\end{proof}

\begin{Lem}
\label{even-lem-em-d}
Let \(\phi\colon M\to N\) be a smooth map and consider a variation of the metric on \(M\) as defined in \eqref{variation-metric-domain}.
Then the following identity holds
\begin{align}
\label{even-lem-em-d-equation}
\int_M&\langle\frac{d}{dt}\big|_{t=0}\Delta^{s-1}\tau(\phi),\Delta^{s-1}\tau(\phi)\rangle\dv \\
\nonumber=&-\sum_{l=1}^{s-1}\int_M\langle\omega,\operatorname{sym}\big(\langle\nabla_{(\cdot)}\Delta^{s-l-1}\tau(\phi),\nabla_{(\cdot)}\Delta^{s+l-2}\tau(\phi)\rangle\big)\rangle\dv \\
\nonumber&-\frac{1}{2}\sum_{l=1}^{s-1}\int_M\langle g,\omega\rangle\langle \Delta^{s-l}\tau(\phi),\Delta^{s+l-2}\tau(\phi)\rangle\dv \\
\nonumber&+\frac{1}{2}\sum_{l=1}^{s-1}\int_M\langle g,\omega\rangle\langle \nabla\Delta^{s-l-1}\tau(\phi),\nabla\Delta^{s+l-2}\tau(\phi)\rangle\dv\\
\nonumber&+\int_M\langle\omega,\operatorname{sym}\big(\langle d\phi(\cdot),\nabla_{(\cdot)}\Delta^{2s-2}\tau(\phi)\rangle\big)\rangle\dv \\
\nonumber&-\frac{1}{2}\int_M\langle g,\omega\rangle\langle\tau(\phi),\Delta^{2s-2}\tau(\phi)\rangle\dv \\
\nonumber&-\frac{1}{2}\int_M\langle g,\omega\rangle\langle d\phi(e_k),\nabla_k\Delta^{2s-2}\tau(\phi)\rangle\dv,
\end{align}
where \(\operatorname{sym}\) denotes symmetrization with respect to the two free slots.
\end{Lem}

\begin{proof}
We insert \eqref{variation-connection-laplacian} into the first term on the right-hand side of \eqref{even-lem-em-c}
and find
\begin{align*}
\int_M\langle&\big(\frac{d}{dt}\big|_{t=0}\Delta\big)\Delta^{s-l-1}\tau(\phi),\Delta^{s+l-2}\tau(\phi)\rangle\dv\\
=&\int_M\omega^{ij}\langle\nabla_i\nabla_j\Delta^{s-l-1}\tau(\phi),\Delta^{s+l-2}\tau(\phi)\rangle\dv 
+\int_M\langle(\nabla^i\omega^{j}_{~i})\nabla_j\Delta^{s-l-1}\tau(\phi),\Delta^{s+l-2}\tau(\phi)\rangle\dv \\
&-\frac{1}{2}\int_M\langle(\nabla^k\tr\omega)\nabla_k\Delta^{s-l-1}\tau(\phi),\Delta^{s+l-2}\tau(\phi)\rangle\dv \\
=&-\int_M\langle\omega,\operatorname{sym}\big(\langle\nabla_{(\cdot)}\Delta^{s-l-1}\tau(\phi),\nabla_{(\cdot)}\Delta^{s+l-2}\tau(\phi)\rangle\big)\rangle\dv \\
&-\frac{1}{2}\int_M\langle g,\omega\rangle\langle \Delta^{s-l}\tau(\phi),\Delta^{s+l-2}\tau(\phi)\rangle\dv 
+\frac{1}{2}\int_M\langle g,\omega\rangle\langle \nabla\Delta^{s-l-1}\tau(\phi),\nabla\Delta^{s+l-2}\tau(\phi)\rangle\dv.
\end{align*}
Moreover, we combine \eqref{variation-tension-field} with the second term on the right hand side of \eqref{even-lem-em-c} yielding
\begin{align*}
\int_M&\langle\frac{d}{dt}\big|_{t=0}\tau(\phi),\Delta^{2s-2}\tau(\phi)\rangle\dv \\
=&-\int_M\omega^{ij}\langle(\nabla d\phi)_{ij},\Delta^{2s-2}\tau(\phi)\rangle\dv 
-\int_M\langle(\nabla_i\omega^{ki})d\phi(e_k),\Delta^{2s-2}\tau(\phi)\rangle\dv \\
&+\frac{1}{2}\int_M\langle(\nabla^k\tr\omega)d\phi(e_k),\Delta^{2s-2}\tau(\phi)\rangle\dv \\
=&\int_M\langle\omega,\operatorname{sym}\big(\langle d\phi(\cdot),\nabla_{(\cdot)}\Delta^{2s-2}\tau(\phi)\rangle\big)\rangle\dv
-\frac{1}{2}\int_M\langle g,\omega\rangle\langle\tau(\phi),\Delta^{2s-2}\tau(\phi)\rangle\dv \\
&-\frac{1}{2}\int_M\langle g,\omega\rangle\langle d\phi(e_k),\nabla_k\Delta^{2s-2}\tau(\phi)\rangle\dv.
\end{align*}
The claim follows by combining both equations.
\end{proof}

After this preliminary work we can define the stress-energy tensor for polyharmonic maps of even order as follows
\begin{align}
\label{even-energy-momentum-tensor}
S_{2s}(X,Y):=&g(X,Y)\bigg(\frac{1}{2}|\Delta^{s-1}\tau(\phi)|^2-\langle\tau(\phi),\Delta^{2s-2}\tau(\phi)\rangle-\langle d\phi,\nabla\Delta^{2s-2}\tau(\phi)\rangle \\
\nonumber&+\sum_{l=1}^{s-1}(-\langle\Delta^{s-l}\tau(\phi),\Delta^{s+l-2}\tau(\phi)\rangle+\langle\nabla\Delta^{s-l-1}\tau(\phi),\nabla\Delta^{s+l-2}\tau(\phi)\rangle)
\bigg) \\
\nonumber&-\sum_{l=1}^{s-1}(\langle\nabla_X\Delta^{s-l-1}\tau(\phi),\nabla_Y\Delta^{s+l-2}\tau(\phi)\rangle+\langle\nabla_Y\Delta^{s-l-1}\tau(\phi),\nabla_X\Delta^{s+l-2}\tau(\phi)\rangle) \\
\nonumber&+\langle d\phi(X),\nabla_Y\Delta^{2s-2}\tau(\phi)\rangle
+\langle d\phi(Y),\nabla_X\Delta^{2s-2}\tau(\phi)\rangle.
\end{align}

\begin{Prop}
Let \(\phi\colon M\to N\) be a smooth map and consider a variation of the metric on \(M\) as defined in \eqref{variation-metric-domain}.
Then the variation of the energy functional \eqref{poly-energy-even} with respect to the metric on the domain is given by
\begin{align}
\frac{d}{dt}\big|_{t=0}\int_M|\Delta^{s-1}\tau(\phi)|^2\dv=\int_M\langle S_{2s},\omega\rangle\dv,
\end{align}
where the symmetric two-tensor \(S_{2s}\) is defined in \eqref{even-energy-momentum-tensor}.
\end{Prop}
\begin{proof}
The claim follows from Lemma \ref{even-lem-em-a}, Lemma \ref{even-lem-em-d} and 
\begin{align*}
\operatorname{sym}\big(\langle d\phi(\cdot),\nabla_{(\cdot)}\Delta^{2s-2}\tau(\phi)\rangle\big)(X,Y)
=\frac{1}{2}\big(\langle d\phi(X),\nabla_{Y}\Delta^{2s-2}\tau(\phi)\rangle+\langle d\phi(Y),\nabla_{X}\Delta^{2s-2}\tau(\phi)\rangle\big).
\end{align*}
The second symmetrization in \eqref{even-lem-em-d-equation} is calculated in the same way.
\end{proof}

It can be directly seen that \(S_{2s}(X,Y)\) is symmetric.
In the following, we will confirm that the stress-energy tensor also satisfies a
conservation law.

\begin{Prop}
Let \(\phi\colon M\to N\) be a smooth map.
Then the stress-energy tensor \eqref{even-energy-momentum-tensor} satisfies the following conservation law
\begin{align}
\operatorname{div}S_{2s}=-\langle\tau_{2s}(\phi),d\phi\rangle.
\end{align}
In particular, the stress-energy tensor is divergence-free whenever \(\phi\) is a solution of \eqref{tension-2s}.
\end{Prop}
\begin{proof}
Throughout the proof we choose a local orthonormal basis \(\{e_i\},i=1,\ldots,m\) of \(TM\) satisfying \(\nabla_{e_i}e_j=0,i,j=1,\ldots,m\) 
around a given point \(p\in M\) and set
\begin{align*}
S_{2s}(e_i,e_j):=&g_{ij}\bigg(\frac{1}{2}|\Delta^{s-1}\tau(\phi)|^2-\langle\tau(\phi),\Delta^{2s-2}\tau(\phi)\rangle-\langle d\phi,\nabla\Delta^{2s-2}\tau(\phi)\rangle \\
\nonumber&+\sum_{l=1}^{s-1}(-\langle\Delta^{s-l}\tau(\phi),\Delta^{s+l-2}\tau(\phi)\rangle+\langle\nabla\Delta^{s-l-1}\tau(\phi),\nabla\Delta^{s+l-2}\tau(\phi)\rangle)
\bigg) \\
\nonumber&-\sum_{l=1}^{s-1}(\langle\nabla_i\Delta^{s-l-1}\tau(\phi),\nabla_j\Delta^{s+l-2}\tau(\phi)\rangle+\langle\nabla_j\Delta^{s-l-1}\tau(\phi),\nabla_i\Delta^{s+l-2}\tau(\phi)\rangle) \\
\nonumber&+\langle d\phi(e_i),\nabla_j\Delta^{2s-2}\tau(\phi)\rangle
+\langle d\phi(e_j),\nabla_i\Delta^{2s-2}\tau(\phi)\rangle.
\end{align*}
By a direct calculation we find
\begin{align*}
\nabla^jS_{2s}(e_i,e_j)=&-\langle d\phi(e_i),\Delta^{2s-1}\tau(\phi)\rangle+\langle d\phi(e_j),(\nabla_j\nabla_i-\nabla_i\nabla_j)\Delta^{2s-2}\tau(\phi)\rangle \\
&+\sum_{l=1}^{s-1}\big(\langle(\nabla_i\nabla_j-\nabla_j\nabla_i)\Delta^{s-l-1}\tau(\phi),\nabla_j\Delta^{s+l-2}\tau(\phi)\rangle \\
&\hspace{0.5cm}+\langle\nabla_j\Delta^{s-l-1}\tau(\phi),(\nabla_i\nabla_j-\nabla_j\nabla_i)\Delta^{s+l-2}\tau(\phi)\rangle\big) \\
&+\langle\nabla_i\Delta^{s-1}\tau(\phi),\Delta^{s-1}\tau(\phi)\rangle
-\langle\nabla_i\tau(\phi),\Delta^{2s-2}\tau(\phi)\rangle \\
&+\sum_{l=1}^{s-1}\big(\langle\nabla_i\Delta^{s-l-1}\tau(\phi),\Delta^{s+l-1}\tau(\phi)\rangle-\langle\Delta^{s+l-2}\tau(\phi),\nabla_i\Delta^{s-l}\tau(\phi)\rangle\big).
\end{align*}
As a next step we rewrite the curvature terms starting with
\begin{align*}
\langle d\phi(e_j),(\nabla_j\nabla_i-\nabla_i\nabla_j)\Delta^{2s-2}\tau(\phi)\rangle=\langle d\phi(e_i),R^N(\Delta^{2s-2}\tau(\phi),d\phi(e_j))d\phi(e_j)\rangle
\end{align*}
and also
\begin{align*}
\langle(\nabla_i\nabla_j-\nabla_j\nabla_i)\Delta^{s-l-1}\tau(\phi)&,\nabla_j\Delta^{s+l-2}\tau(\phi)\rangle \\
&=-\langle R^N(\Delta^{s-l-1}\tau(\phi),\nabla_j\Delta^{s+l-2}\tau(\phi))d\phi(e_j),d\phi(e_i)\rangle, \\
\langle\nabla_j\Delta^{s-l-1}\tau(\phi)&,(\nabla_i\nabla_j-\nabla_j\nabla_i)\Delta^{s+l-2}\tau(\phi)\rangle  \\
&=-\langle R^N(\Delta^{s+l-2}\tau(\phi),\nabla_j\Delta^{s-l-1}\tau(\phi))d\phi(e_j),d\phi(e_i)\rangle.
\end{align*}
Finally, we have the following kind of telescope sum identity
\begin{align*}
\sum_{l=1}^{s-1}(\langle\nabla_i\Delta^{s-l-1}\tau(\phi),&\Delta^{s+l-1}\tau(\phi)\rangle-\langle\nabla_i\Delta^{s-l}\tau(\phi),\Delta^{s+l-2}\tau(\phi)\rangle) \\
&=\langle\nabla_i\tau(\phi),\Delta^{2s-2}\tau(\phi)\rangle
-\langle\nabla_i\Delta^{s-1}\tau(\phi),\Delta^{s-1}\tau(\phi)\rangle.
\end{align*}
The claim follows by combining all equations.
\end{proof}

\subsection{The odd case}
In this section we derive the stress-energy tensor in the case of a polyharmonic map of odd order.
The calculation is similar to the even case and we do not give as many details as before.

\begin{Lem}
\label{odd-lem-em-a}
Let \(\phi\colon M\to N\) be a smooth map and consider a variation of the metric on \(M\) as defined in \eqref{variation-metric-domain}.
Then the following identity holds
\begin{align}
\nonumber\frac{d}{dt}\big|_{t=0}E_{2s+1}(\phi,g_t)=&\frac{1}{2}\int_M|\nabla\Delta^{s-1}\tau(\phi)|^2\langle g,\omega\rangle \dv
-\int_M\omega^{ij}\langle\nabla_i\Delta^{s-1}\tau(\phi),\nabla_j\Delta^{s-1}\tau(\phi)\rangle\dv \\
&+2\int_M\langle\frac{d}{dt}\big|_{t=0}\nabla\Delta^{s-1}\tau(\phi),\nabla\Delta^{s-1}\tau(\phi)\rangle\dv.
\end{align}
\end{Lem}

\begin{proof}
This follows by a direct calculation.
\end{proof}

\begin{Lem}
Let \(\phi\colon M\to N\) be a smooth map and consider a variation of the metric on \(M\) as defined in \eqref{variation-metric-domain}.
Then the following identity holds
\begin{align*}
\nonumber\int_M\langle\frac{d}{dt}\big|_{t=0}\nabla\Delta^{s-1}\tau(\phi),\nabla\Delta^{s-1}\tau(\phi)\rangle\dv=&
\sum_{l=1}^{s-1}\int_M\langle\big(\frac{d}{dt}\big|_{t=0}\Delta\big)\Delta^{s-l-1}\tau(\phi),\Delta^{s+l-1}\tau(\phi)\rangle\dv \\
&+\int_M\langle\frac{d}{dt}\big|_{t=0}\tau(\phi),\Delta^{2s-1}\tau(\phi)\rangle\dv.
\end{align*}
\end{Lem}

\begin{proof}
We note that
\begin{align*}
\frac{d}{dt}\big|_{t=0}\big(\nabla\Delta^{s-1}\tau(\phi)\big)=&\nabla\sum_{l=1}^{s-1}\Delta^{l-1}\big(\frac{d}{dt}\big|_{t=0}\Delta\big)\Delta^{s-l-1}\tau(\phi)
+\nabla\Delta^{s-1}\frac{d}{dt}\big|_{t=0}\tau(\phi) 
\end{align*}
and the claim follows using integration by parts and the self-adjointness of the Laplacian.
\end{proof}

\begin{Lem}
\label{odd-lem-em-d}
Let \(\phi\colon M\to N\) be a smooth map and consider a variation of the metric on \(M\) as defined in \eqref{variation-metric-domain}.
Then the following identity holds
\begin{align}
\int_M&\langle\frac{d}{dt}\big|_{t=0}\nabla\Delta^{s-1}\tau(\phi),\nabla\Delta^{s-1}\tau(\phi)\rangle\dv \\
\nonumber=&-\sum_{l=1}^{s-1}\int_M\langle\omega,\operatorname{sym}\big(\langle\nabla_{(\cdot)}\Delta^{s-l-1}\tau(\phi),\nabla_{(\cdot)}\Delta^{s+l-1}\tau(\phi)\rangle\big)\rangle\dv \\
\nonumber&-\frac{1}{2}\sum_{l=1}^{s-1}\int_M\langle g,\omega\rangle\langle \Delta^{s-l}\tau(\phi),\Delta^{s+l-1}\tau(\phi)\rangle\dv \\
\nonumber&+\frac{1}{2}\sum_{l=1}^{s-1}\int_M\langle g,\omega\rangle\langle \nabla\Delta^{s-l-1}\tau(\phi),\nabla\Delta^{s+l-1}\tau(\phi)\rangle\dv\\
\nonumber&+\int_M\langle\omega,\operatorname{sym}\big(\langle d\phi(\cdot),\nabla_{(\cdot)}\Delta^{2s-1}\tau(\phi)\rangle\big)\rangle\dv \\
\nonumber&-\frac{1}{2}\int_M\langle g,\omega\rangle\langle\tau(\phi),\Delta^{2s-1}\tau(\phi)\rangle\dv \\
\nonumber&-\frac{1}{2}\int_M\langle g,\omega\rangle\langle d\phi(e_k),\nabla_k\Delta^{2s-1}\tau(\phi)\rangle\dv,
\end{align}
where \(\operatorname{sym}\) denotes symmetrization with respect to the two free slots.
\end{Lem}

\begin{proof}
The claim follows similar to the even case, see the proof of Lemma \ref{even-lem-em-d}.
\end{proof}

Again, we are ready to define the stress-energy tensor for polyharmonic maps of odd order by 
\begin{align}
\label{odd-energy-momentum-tensor}
S_{2s+1}(X,Y):=&g(X,Y)\bigg(\frac{1}{2}|\nabla\Delta^{s-1}\tau(\phi)|^2-\langle\tau(\phi),\Delta^{2s-1}\tau(\phi)\rangle-\langle d\phi,\nabla\Delta^{2s-1}\tau(\phi)\rangle \\
\nonumber&+\sum_{l=1}^{s-1}(-\langle\Delta^{s-l}\tau(\phi),\Delta^{s+l-1}\tau(\phi)\rangle+\langle\nabla\Delta^{s-l-1}\tau(\phi),\nabla\Delta^{s+l-1}\tau(\phi)\rangle)
\bigg) \\
\nonumber&-\sum_{l=1}^{s-1}(\langle\nabla_X\Delta^{s-l-1}\tau(\phi),\nabla_Y\Delta^{s+l-1}\tau(\phi)\rangle+\langle\nabla_Y\Delta^{s-l-1}\tau(\phi),\nabla_X\Delta^{s+l-1}\tau(\phi)\rangle) \\
\nonumber&+\langle d\phi(X),\nabla_Y\Delta^{2s-1}\tau(\phi)\rangle
+\langle d\phi(Y),\nabla_X\Delta^{2s-1}\tau(\phi)\rangle\\
\nonumber&-\langle\nabla_X\Delta^{s-1}\tau(\phi),\nabla_Y\Delta^{s-1}\tau(\phi)\rangle.
\end{align}

\begin{Bem}
Note that the structure of the stress-energy tensors \eqref{even-energy-momentum-tensor} and \eqref{odd-energy-momentum-tensor} is 
similar. In both cases the first term arises from the variation of the volume element.
Except from the last contribution in \eqref{odd-energy-momentum-tensor} all terms have the same structure and only differ 
in the order of differentiation. The last term in \eqref{odd-energy-momentum-tensor} arises since we have to vary
the metric one additional time in \eqref{poly-energy-odd}.
\end{Bem}

\begin{Prop}
Let \(\phi\colon M\to N\) be a smooth map and consider a variation of the metric on \(M\) as defined in \eqref{variation-metric-domain}.
Then the variation of the energy functional \eqref{poly-energy-odd} with respect to the metric on the domain is given by
\begin{align}
\frac{d}{dt}\big|_{t=0}\int_M|\nabla\Delta^{s-1}\tau(\phi)|^2\dv=\int_M\langle S_{2s+1},\omega\rangle\dv,
\end{align}
where the symmetric two-tensor \(S_{2s+1}\) is defined in \eqref{odd-energy-momentum-tensor}.
\end{Prop}
\begin{proof}
This follows from Lemma \ref{odd-lem-em-a}, Lemma \ref{odd-lem-em-d} and the same arguments that were used in the even case.
\end{proof}

\begin{Prop}
\label{conservation-stress-energy-odd}
Let \(\phi\colon M\to N\) be a smooth map.
Then the stress-energy tensor \eqref{odd-energy-momentum-tensor} satisfies the following conservation law
\begin{align}
\operatorname{div}S_{2s+1}=-\langle\tau_{2s+1}(\phi),d\phi\rangle.
\end{align}
In particular, the stress-energy tensor is divergence-free whenever \(\phi\) is a solution of \eqref{tension-2s+1}.
\end{Prop}
\begin{proof}
Throughout the proof we choose a local orthonormal basis \(\{e_i\},i=1,\ldots,m\) of \(TM\) satisfying \(\nabla_{e_i}e_j=0,i,j=1,\ldots,m\) 
around a given point \(p\in M\) and set
\begin{align*}
S_{2s+1}(e_i,e_j):=&g_{ij}\bigg(\frac{1}{2}|\nabla\Delta^{s-1}\tau(\phi)|^2-\langle\tau(\phi),\Delta^{2s-1}\tau(\phi)\rangle-\langle d\phi,\nabla\Delta^{2s-1}\tau(\phi)\rangle \\
&+\sum_{l=1}^{s-1}(-\langle\Delta^{s-l}\tau(\phi),\Delta^{s+l-1}\tau(\phi)\rangle+\langle\nabla\Delta^{s-l-1}\tau(\phi),\nabla\Delta^{s+l-1}\tau(\phi)\rangle)
\bigg) \\
&-\sum_{l=1}^{s-1}(\langle\nabla_i\Delta^{s-l-1}\tau(\phi),\nabla_j\Delta^{s+l-1}\tau(\phi)\rangle+\langle\nabla_j\Delta^{s-l-1}\tau(\phi),\nabla_i\Delta^{s+l-1}\tau(\phi)\rangle) \\
&+\langle d\phi(e_i),\nabla_j\Delta^{2s-1}\tau(\phi)\rangle
+\langle d\phi(e_j),\nabla_i\Delta^{2s-1}\tau(\phi)\rangle\\
&-\langle\nabla_i\Delta^{s-1}\tau(\phi),\nabla_j\Delta^{s-1}\tau(\phi)\rangle.
\end{align*}
Now, a direct calculation yields
\begin{align*}
\nabla^jS_{2s+1}(e_i,e_j)=&-\langle d\phi(e_i),\Delta^{2s}\tau(\phi)\rangle
+\langle d\phi(e_j),(\nabla_j\nabla_i-\nabla_i\nabla_j)\Delta^{2s-1}\tau(\phi)\rangle\\
&+\langle\nabla_j\Delta^{s-1}\tau(\phi),(\nabla_i\nabla_j-\nabla_j\nabla_i)\Delta^{s-1}\tau(\phi)\rangle \\
&+\sum_{l=1}^{s-1}\big(
\langle\nabla_j\Delta^{s+l-1}\tau(\phi),(\nabla_i\nabla_j-\nabla_j\nabla_i)\Delta^{s-l-1}\tau(\phi)\rangle \\
&\hspace{0.5cm}+\langle\nabla_j\Delta^{s-l-1}\tau(\phi),(\nabla_i\nabla_j-\nabla_j\nabla_i)\Delta^{s+l-1}\tau(\phi)\rangle
\big)\\
&+\sum_{l=1}^{s-1}\big(-\langle\Delta^{s+l-1}\tau(\phi),\nabla_i\Delta^{s-l}\tau(\phi)\rangle
+\langle\Delta^{s+l}\tau(\phi),\nabla_i\Delta^{s-l-1}\tau(\phi)\rangle\big) \\
&-\langle\nabla_i\tau(\phi),\Delta^{2s-1}\tau(\phi)\rangle+\langle\Delta^s\tau(\phi),\nabla_i\Delta^{s-1}\tau(\phi)\rangle.
\end{align*}
As in the even case we have the following ``telescope sum'' identity
\begin{align*}
\sum_{l=1}^{s-1}\big(\langle\Delta^{s+l}\tau(\phi),\nabla_i\Delta^{s-l-1}\tau(\phi)\rangle&
-\langle\Delta^{s+l-1}\tau(\phi),\nabla_i\Delta^{s-l}\tau(\phi)\rangle\big) \\
&=\langle\nabla_i\tau(\phi),\Delta^{2s-1}\tau(\phi)\rangle-\langle\Delta^s\tau(\phi),\nabla_i\Delta^{s-1}\tau(\phi)\rangle.
\end{align*}
Again, we have to rewrite the curvature terms as follows
\begin{align*}
\langle d\phi(e_j),(\nabla_j\nabla_i-\nabla_i\nabla_j)\Delta^{2s-1}\tau(\phi)\rangle
=\langle R^N(\Delta^{2s-1}\tau(\phi),d\phi(e_j))d\phi(e_j),d\phi(e_i)\rangle.
\end{align*}
Moreover, we have
\begin{align*}
\langle\nabla_j\Delta^{s-1}\tau(\phi),(\nabla_i\nabla_j-\nabla_j\nabla_i)\Delta^{s-1}\tau(\phi)\rangle
=-\langle R^N(\Delta^{s-1}\tau(\phi),\nabla_j\Delta^{s-1}\tau(\phi))d\phi(e_j),d\phi(e_i)\rangle
\end{align*}
and also
\begin{align*}
\langle\nabla_j\Delta^{s-l-1}\tau(\phi)&,(\nabla_i\nabla_j-\nabla_j\nabla_i)\Delta^{s+l-1}\tau(\phi)\rangle \\
&=-\langle R^N(\Delta^{s+l-1}\tau(\phi),\nabla_j\Delta^{s-l-1}\tau(\phi))d\phi(e_j),d\phi(e_i)\rangle,\\
\langle\nabla_j\Delta^{s+l-1}\tau(\phi)&,(\nabla_i\nabla_j-\nabla_j\nabla_i)\Delta^{s-l-1}\tau(\phi)\rangle \\
&=-\langle R^N(\Delta^{s-l-1}\tau(\phi),\nabla_j\Delta^{s+l-1}\tau(\phi))d\phi(e_j),d\phi(e_i)\rangle.
\end{align*}
The claim follows by combining all equations.
\end{proof}

\begin{Bem}
Although the Euler-Lagrange equations for polyharmonic maps \eqref{tension-2s} and \eqref{tension-2s+1}
contain various curvature contributions from the target the stress-energy tensors \eqref{even-energy-momentum-tensor}
and \eqref{odd-energy-momentum-tensor} do not.
\end{Bem}

\subsection{Invariance of the energy functional}
In the previous section we have seen that the stress-energy tensors \eqref{even-energy-momentum-tensor}
and \eqref{odd-energy-momentum-tensor} are conserved whenever we have a polyharmonic map,
that is a solution of \eqref{tension-2s} or \eqref{tension-2s+1}.
The existence of these conservation laws does not come as a surprise,
instead it is a direct consequence of Noether's theorem and the invariance 
of the energy functionals \eqref{poly-energy-even} and \eqref{poly-energy-odd} under diffeomorphisms on the domain.
We will make this statement more precise in this section.

First of all, we will recall the invariance of the tension field under diffeomorphisms on the domain.
In order to emphasize the dependence of the tension field on the domain metric \(g\)
we write \(\tau_g(\phi)\).

\begin{Lem}
\label{lemma-invariance-tension-field}
Let \(u\colon M\to M\) be a diffeomorphism. 
Then we have 
\begin{align}
\label{invariance-tension-field}
\tau_g(\phi)=\tau_{u^\ast g}(\phi\circ u).
\end{align}
\end{Lem}
\begin{proof}
Again, we will prove this statement making use of a local calculation.
To this end let \((U,x^i)\) be a local coordinate chart on \(M\) and \((V,y^\alpha)\) be a local coordinate chart on \(N\)
such that \(\phi(U)\subset V\).
Then we find
\begin{align*}
\tau^\alpha_{u^\ast g}(\phi\circ u)=&(u^\ast g)^{ij}\big(
\frac{\partial^2(\phi\circ u)^\alpha}{\partial x^i\partial x^j}-\Gamma^k_{ij}(u^\ast g)\frac{\partial(\phi\circ u)^\alpha}{\partial x^k}
+\Gamma^\alpha_{\beta\gamma}\frac{\partial(\phi\circ u)^\beta}{\partial x^i}\frac{\partial(\phi\circ u)^\gamma}{\partial x^j}
\big).
\end{align*}
By a direct calculation we obtain the identities
\begin{align*}
\frac{\partial(\phi\circ u)^\alpha}{\partial x^i}=&\frac{\partial\phi^\alpha}{\partial u^r}\frac{\partial u^r}{\partial x^i},\\
\frac{\partial^2(\phi\circ u)^\alpha}{\partial x^i\partial x^j}=&\frac{\partial^2\phi^\alpha}{\partial u^r\partial u^s}\frac{\partial u^r}{\partial x^i}
\frac{\partial u^s}{\partial x^j}+\frac{\partial\phi^\alpha}{\partial u^r}\frac{\partial^2u^r}{\partial x^i\partial x^j}.
\end{align*}
In addition, it is well known that the Christoffel symbols satisfy the following formula
\begin{align}
\label{trafo-christoffel}
\Gamma^k_{ij}(u^\ast g)=\frac{\partial u^r}{\partial x^i}\frac{\partial u^s}{\partial x^j}\frac{\partial x^k}{\partial u^t}\Gamma^t_{rs}(g)+\frac{\partial^2 u^t}{\partial x^i\partial x^j}\frac{\partial x^k}{\partial u^t}.
\end{align}
Applying these identities we find
\begin{align*}
\tau^\alpha_{u^\ast g}(\phi\circ u)&=
g^{ab}\frac{\partial x^i}{\partial u^a}\frac{\partial x^j}{\partial u^b}
\big(
\frac{\partial^2\phi^\alpha}{\partial u^r\partial u^s}\frac{\partial u^r}{\partial x^i}
\frac{\partial u^s}{\partial x^j}+\frac{\partial\phi^\alpha}{\partial u^r}\frac{\partial^2u^r}{\partial x^i\partial x^j} 
-\frac{\partial u^r}{\partial x^i}\frac{\partial u^s}{\partial x^j}\frac{\partial x^k}{\partial u^t}\Gamma^t_{rs}(g)\frac{\partial\phi^\alpha}{\partial u^p}\frac{\partial u^p}{\partial x^k} \\
&\hspace{2.5cm}-\frac{\partial^2 u^t}{\partial x^i\partial x^j}\frac{\partial x^k}{\partial u^t}\frac{\partial\phi^\alpha}{\partial u^p}\frac{\partial u^p}{\partial x^k}
+\Gamma^\alpha_{\beta\gamma}\frac{\partial\phi^\beta}{\partial u^r}\frac{\partial\phi^\gamma}{\partial u^s}\frac{\partial u^r}{\partial x^i}\frac{\partial u^s}{\partial x^j}
\big) \\
&=g^{rs}\big(
\frac{\partial^2\phi^\alpha}{\partial u^r\partial u^s}
-\Gamma^t_{rs}(g)\frac{\partial\phi^\alpha}{\partial u^t}
+\Gamma^\alpha_{\beta\gamma}\frac{\partial\phi^\beta}{\partial u^r}\frac{\partial\phi^\gamma}{\partial u^s}
\big)\\
&=\tau^\alpha_g(\phi),
\end{align*}
which completes the proof.
\end{proof}

\begin{Bem}
It is to be expected that \eqref{invariance-tension-field} holds as 
the Dirichlet energy is invariant under diffeomorphisms \(u\colon M\to M\) in the following sense
\begin{align*}
E(\phi\circ u,u^\ast g)=E(\phi,g).
\end{align*}
Hence, one should expect that the critical points of this functional share the same invariance.
\end{Bem}

In order to highlight the dependence of the connection Laplacian on the vector bundle \(\phi^\ast TN\)
on both metric of the domain \(g\) and map \(\phi\) we write \(\Delta^\phi_g\).

\begin{Lem}
Let \(u\colon M\to M\) be a diffeomorphism. 
Then we have 
\begin{align}
\label{invariance-connection-laplacian}
\Delta^\phi_g=\Delta^{\phi\circ u}_{u\ast g}.
\end{align}
\end{Lem}

\begin{proof}
As in the proof of Lemma \ref{lemma-invariance-tension-field} we perform a calculation
in local coordinates to prove the statement.
Recall that in terms of local coordinates we have
\begin{align*}
\Delta^\phi_{g}=-g^{ij}(\nabla^\phi_i\nabla^\phi_j-\Gamma^k_{ij}(g)\nabla^\phi_k).
\end{align*}
In order to show the invariance of the connection Laplacian under diffeomorphisms we calculate
\begin{align*}
(u^\ast g)^{ij}\nabla_i^{\phi\circ u}\nabla_j^{\phi\circ u}=&
(u^\ast g)^{ij}\nabla_i^{\phi\circ u}\big(\frac{\partial\phi^\alpha}{\partial u^r}\frac{\partial u^r}{\partial x^j}\nabla^N_{\partial y^\alpha}\big)\\
=&(u^\ast g)^{ij}\nabla_i^{\phi\circ u}\big(\frac{\partial u^r}{\partial x^j}\nabla^\phi_{r}\big) \\
=&(u^\ast g)^{ij}\frac{\partial^2u^r}{\partial x^i\partial x^j}\nabla^\phi_r
+(u^\ast g)^{ij}\frac{\partial u^r}{\partial x^j}\frac{\partial u^s}{\partial x^i}\nabla^\phi_{r}\nabla^\phi_{s}\\
=&(u^\ast g)^{ij}\frac{\partial^2u^r}{\partial x^i\partial x^j}\nabla^\phi_r
+g^{rs}\nabla^\phi_{r}\nabla^\phi_{s}.
\end{align*}
Here, we used \(\nabla^N\) to denote the Levi-Civita connection on \(N\).
In addition, we again apply the formula for the transformation of the Christoffel symbols \eqref{trafo-christoffel}
and find
\begin{align*}
\Gamma^k_{ij}(u^\ast g)\nabla_k^{\phi\circ u}=&
\frac{\partial u^r}{\partial x^i}\frac{\partial u^s}{\partial x^j}\frac{\partial x^k}{\partial u^t}\Gamma^t_{rs}(g)\frac{\partial\phi^\alpha}{\partial u^q}\frac{\partial u^q}{\partial x^k}\nabla^N_{\partial y_\alpha}
+\frac{\partial^2 u^t}{\partial x^i\partial x^j}\frac{\partial x^k}{\partial u^t}\frac{\partial\phi^\alpha}{\partial u^q}\frac{\partial u^q}{\partial x^k}\nabla^N_{\partial y_\alpha}\\
=&\frac{\partial u^r}{\partial x^i}\frac{\partial u^s}{\partial x^j}\Gamma^t_{rs}(g)\nabla_t^\phi
+\frac{\partial^2 u^r}{\partial x^i\partial x^j}\nabla_r^\phi.
\end{align*}
The claim follows by combining both equations.
\end{proof}

Summarizing the considerations in this section we can draw the following conclusion:
\begin{Satz}
Let \(u\colon M\to M\) be a diffeomorphism.
Then the energy functionals \eqref{poly-energy-even} and \eqref{poly-energy-odd} are 
invariant under diffeomorphisms on the domain in the following sense
\begin{align*}
E_{2s}(\phi\circ u,u^\ast g)=E_{2s}(\phi,g),\qquad E_{2s+1}(\phi\circ u,u^\ast g)=E_{2s+1}(\phi,g),
\end{align*}
where \(s=1,2,\ldots\).
\end{Satz}
\begin{proof}
In the even case this follows directly from \eqref{invariance-tension-field} and \eqref{invariance-connection-laplacian}.
In the odd case we note that
\begin{align*}
E_{2s+1}(\phi\circ u,u^\ast g)
&=\int_M(u^\ast g)^{ij}\langle\nabla^{\phi\circ u}_i(\Delta^{\phi\circ u}_{u^\ast g}\big)^{s-1}\tau_{u^\ast g}(u\circ\phi),
\nabla^{\phi\circ u}_j(\Delta^{\phi\circ u}_{u^\ast g}\big)^{s-1}\tau_{u^\ast g}(u\circ\phi)\rangle\dv_{u^\ast g} \\
&=\int_M(u^\ast g)^{ij}\langle\nabla^{\phi\circ u}_i(\Delta^\phi_{g}\big)^{s-1}\tau_{g}(\phi),
\nabla^{\phi\circ u}_j(\Delta^{\phi}_{g}\big)^{s-1}\tau_{g}(\phi)\rangle\dv \\
&=\int_Mg^{ij}\langle\nabla^{\phi}_i(\Delta^\phi_{g}\big)^{s-1}\tau_{g}(\phi),
\nabla^{\phi}_j(\Delta^{\phi}_{g}\big)^{s-1}\tau_{g}(\phi)\rangle\dv \\
&=E_{2s+1}(\phi,g),
\end{align*}
where we used \eqref{invariance-tension-field} and \eqref{invariance-connection-laplacian} in the first step.
The proof is now complete.
\end{proof}

\section{Maps with vanishing stress-energy tensor}
In this section we apply the stress energy-tensor for polyharmonic maps,
that is \eqref{even-energy-momentum-tensor} in the even case and \eqref{odd-energy-momentum-tensor} in the odd case,
in order to characterize the qualitative behavior of polyharmonic maps.

\subsection{Liouville-type results from vanishing stress-energy tensor: The even case}
We find that the trace of \eqref{even-energy-momentum-tensor} is given by
\begin{align*}
\tr S_{2s}=&\frac{m}{2}|\Delta^{s-1}\tau(\phi)|^2+(2-m)\langle d\phi,\nabla\Delta^{2s-2}\tau(\phi)\rangle
-m\langle\tau(\phi),\Delta^{2s-2}\tau(\phi)\rangle \\
&-m\sum_{l=1}^{s-1}
\langle\Delta^{s-l}\tau(\phi),\Delta^{s+l-2}\tau(\phi)\rangle
+(m-2)\sum_{l=1}^{s-1}\langle\nabla\Delta^{s-l-1}\tau(\phi),\nabla\Delta^{s+l-2}\tau(\phi)\rangle.
\end{align*}

\begin{Satz}
Suppose that \((M,g)\) is a closed Riemannian manifold and let \(\phi\colon M\to N\) be a smooth map with \(S_{2s}=0\).
If \(m\neq 4s\) then \(\phi\) must be harmonic.
\end{Satz}
\begin{proof}
By assumption the stress-energy tensor \eqref{even-energy-momentum-tensor} vanishes.
Consequently, it is also trace-free and after using integration by parts several times 
we find
\begin{align*}
0=\int_M\tr S_{2s}\dv=(\frac{m}{2}-2s)\int_M|\Delta^{s-1}\tau(\phi)|^2\dv.
\end{align*}
Again, by assumption \(m\neq 4s\) which implies that \(\Delta^{s-1}\tau(\phi)=0\).
Using integration by parts once more, we find
\begin{align*}
0=\int_M\langle \Delta^{s-1}\tau(\phi),\Delta^{s-2}\tau(\phi)\rangle\dv=\int_M|\nabla\Delta^{s-2}\tau(\phi)|^2\dv
\end{align*}
yielding that \(\nabla\Delta^{s-2}\tau(\phi)=0\). 

Then, we calculate
\begin{align*}
0=\int_M\langle\nabla\Delta^{s-2}\tau(\phi),\nabla\Delta^{s-3}\tau(\phi)\rangle\dv=\int_M|\Delta^{s-2}\tau(\phi)|^2\dv.
\end{align*}
By iterating this procedure we finally obtain that \(\tau(\phi)=0\) proving the claim.
\end{proof}

\subsection{Liouville-type results from vanishing stress-energy tensor: the odd case}
Now, we find that the trace of \eqref{odd-energy-momentum-tensor} is given by
\begin{align*}
\tr S_{2s+1}=&(\frac{m}{2}-1)|\nabla\Delta^{s-1}\tau(\phi)|^2+(2-m)\langle d\phi,\nabla\Delta^{2s-1}\tau(\phi)\rangle
-m\langle\tau(\phi),\Delta^{2s-1}\tau(\phi)\rangle \\
&-m\sum_{l=1}^{s-1}\langle\Delta^{s-l}\tau(\phi),\Delta^{s+l-1}\tau(\phi)\rangle
+(m-2)\sum_{l=1}^{s-1}\langle\nabla\Delta^{s-l-1}\tau(\phi),\nabla\Delta^{s+l-1}\tau(\phi)\rangle.
\end{align*}

\begin{Satz}
Suppose that \((M,g)\) is a closed Riemannian manifold and let \(\phi\colon M\to N\) be a smooth map with \(S_{2s+1}=0\).
If \(m\neq 4s+2\) then \(\phi\) must be harmonic.
\end{Satz}
\begin{proof}
By assumption the stress-energy tensor \eqref{odd-energy-momentum-tensor} vanishes.
Consequently, it is also trace-free and after using integration by parts several times 
we find
\begin{align*}
0=\int_M\tr S_{2s+1}\dv=(\frac{m}{2}-2s-1)\int_M|\nabla\Delta^{s-1}\tau(\phi)|^2\dv.
\end{align*}
By assumption \(m\neq 4s+2\) which implies that \(\nabla\Delta^{s-1}\tau(\phi)=0\).
Using integration by parts once more, we find
\begin{align*}
0=\int_M\langle\nabla\Delta^{s-1}\tau(\phi),\nabla\Delta^{s-2}\tau(\phi)\rangle\dv=\int_M|\Delta^{s-1}\tau(\phi)|^2\dv
\end{align*}
yielding that \(\Delta^{s-1}\tau(\phi)=0\). By iterating this procedure we finally obtain that \(\tau(\phi)=0\)
proving the claim.
\end{proof}

\section{A classification result for triharmonic maps}
In this section we apply the stress-energy tensor to obtain a classification result for \emph{triharmonic maps}.
These are critical points of the trienergy, which is given by
\begin{align}
\label{tri-energy}
E(\phi)=\int_M|\nabla\tau(\phi)|^2\dv
\end{align}
and corresponds to \eqref{poly-energy-odd} with \(s=1\).

A critical point of \eqref{tri-energy} is called \emph{triharmonic map} and solves the following
partial differential equation of sixth order
\begin{align}
\label{triharmonic-map}
\Delta^2\tau(\phi)=R^N(\Delta\tau(\phi),d\phi(e_i))d\phi(e_i)+R^N(\nabla_i\tau(\phi),\tau(\phi))d\phi(e_i).
\end{align}

The stress-energy tensor obtained from the trienergy \eqref{tri-energy} is given by
\begin{align}
\label{stress-energy-triharmonic}
S_3(e_i,e_j)=&g_{ij}(\frac{1}{2}|\nabla\tau(\phi)|^2-\langle\tau(\phi),\Delta\tau(\phi)\rangle-\langle d\phi,\nabla\Delta\tau(\phi)\rangle) \\
\nonumber&-\langle\nabla_i\tau(\phi),\nabla_j\tau(\phi)\rangle+\langle d\phi(e_i),\nabla_j\Delta\tau(\phi)\rangle
+\langle d\phi(e_j),\nabla_i\Delta\tau(\phi)\rangle.
\end{align}

Currently, there are only few results available in the literature on triharmonic maps.
Explicit solutions to the triharmonic map equation have been considered in \cite{MR3403738,MR3711937,MR3790367}.
Triharmonic immersions into Riemannian manifolds of non-positive curvature have been studied in \cite{MR3371364}.
A monotonicity formula for extrinsic triharmonic maps could be achieved in \cite{MR3624536}.
A triharmonic geometric heat flow for surfaces in \(\R^3\) has been investigated in \cite{MR3672995}.

In the following we will show that triharmonic maps of finite energy from Euclidean space must be harmonic or trivial.
This extends a previous non-existence result for biharmonic maps of finite energy, see \cite[Theorem 3.4]{MR2604617}.
However, here, we are using a different method of proof.

\begin{Satz}
\label{liouville-triharmonic}
Let \(\phi\colon\R^m\to N\) be a smooth triharmonic map and \(m\neq 6\).
Suppose that 
\begin{align}
\label{thm-finiteness-triharmonic}
\int_{\R^m}(|d\phi|^2+|\nabla d\phi|^2+|\nabla^2d\phi|^2)\dv<\infty.
\end{align}
If \(m=2\) then \(\phi\) must be harmonic, if \(m>2\) then \(\phi\) must be constant.
\end{Satz}

As in the case of biharmonic maps we have to impose the finiteness of
the full third covariant derivative of \(\phi\) which is a stronger assumption
as demanding finite trienergy. Note that \(|\nabla\tau(\phi)|^2\leq m|\nabla^2d\phi|^2\).

\begin{proof}
For \(R>0\) let \(\eta\in C_0^\infty(\mathbb{R})\) be a smooth cut-off function satisfying \(\eta=1\) for \(|z|\leq R\),
\(\eta=0\) for \(|z|\geq 2R\) and \(|\eta^l(z)|\leq\frac{C}{R^l},l=1,2,3\). In addition, 
we set \(Y(x):=x\eta(r)\in C_0^\infty(\mathbb{R}^m,\mathbb{R}^m)\) with \(r=|x|\).
Hence, we find
\[
\frac{\partial Y_i}{\partial x^j}=\delta_{ij}\eta(r)+\frac{x_i x_j}{r}\eta'(r).
\]
Due to Proposition \ref{conservation-stress-energy-odd} we know that the stress-energy tensor \eqref{stress-energy-triharmonic} 
is divergence-free whenever we assume that \(\phi\) is a solution of \eqref{triharmonic-map}.
Hence, we find
\begin{align*}
0=-\int_{\R^m}\langle Y,\operatorname{div}S_3\rangle\dv=\int_{\R^m}\frac{\partial Y_i}{\partial x^j}S_3(e_i,e_j)\dv.
\end{align*}
By a direct computation we find
\begin{align*}
\int_{\R^m} S_3(e_i,e_j)\delta_{ij}\eta(r)\dv=& \int_{\R^m}\eta(r)\big((\frac{m}{2}-1)|\nabla\tau(\phi)|^2 \\
&+(2-m)\langle\nabla\Delta\tau(\phi),d\phi\rangle-m\langle\Delta\tau(\phi),\tau(\phi)\rangle\big)\dv.
\end{align*}

Now, we calculate
\begin{align*}
\int_{\R^m}\eta(r)\langle\Delta\tau(\phi),\tau(\phi)\rangle\dv=\int_{\R^m}\eta(r)|\nabla\tau(\phi)|^2\dv+\int_{\R^m}\eta'(r)\frac{x_j}{r}\langle\nabla_j\tau(\phi),\tau(\phi)\rangle\dv
\end{align*}
and also
\begin{align*}
\int_{\R^m}\eta(r)\langle\nabla\Delta\tau(\phi),d\phi\rangle\dv=&-\int_{\R^m}\eta(r)|\nabla\tau(\phi)|^2\dv-\int_{\R^m}\eta'(r)\frac{x_j}{r}\langle\nabla_j\tau(\phi),\tau(\phi)\rangle\dv \\
&-\int_{\R^m}\eta'(r)\frac{x_j}{r}\langle\Delta\tau(\phi),d\phi(e_j)\rangle\dv.
\end{align*}

This allows us to deduce
\begin{align*}
\int_{\R^m}S_3(e_i,e_j)\delta_{ij}\eta(r)\dv=& \int_{\R^m}\eta(r)(\frac{m}{2}-3)|\nabla\tau(\phi)|^2\dv
-2\int_{\R^m}\eta'(r)\frac{x_j}{r}\langle\nabla_j\tau(\phi),\tau(\phi)\rangle\dv \\
&-(2-m)\int_{\R^m}\eta'(r)\frac{x_j}{r}\langle\Delta\tau(\phi),d\phi(e_j)\rangle\dv.
\end{align*}

In addition, we find
\begin{align*}
\int_{\R^m}S_3(e_i,e_j)\frac{x_i x_j}{r}\eta'(r)\dv=\int_{\R^m}&\eta'(r)\big(\frac{1}{2}r|\nabla\tau(\phi)|^2-r\langle\Delta\tau(\phi),\tau(\phi)\rangle 
-r\langle\nabla\Delta\tau(\phi),d\phi\rangle \\
&-\frac{x_ix_j}{r}\langle\nabla_i\tau(\phi),\nabla_j\tau(\phi)\rangle+2\frac{x_ix_j}{r}\langle\nabla_i\Delta\tau(\phi),d\phi(e_j)\rangle 
\big)\dv.
\end{align*}

Note that
\begin{align*}
\int_{\R^m}\eta'(r)r\big(\langle\nabla\Delta\tau(\phi),d\phi\rangle+\langle\Delta\tau(\phi),\tau(\phi)\rangle\big)\dv
=-\int_{\R^m}\big(\eta''(r)x_j+\eta'(r)\frac{x_j}{r}\big)\langle\Delta\tau(\phi),d\phi(e_j)\rangle\dv
\end{align*}
and also
\begin{align*}
\int_{\R^m}\eta'(r)\frac{x_ix_j}{r}\langle\nabla_i\Delta\tau(\phi),d\phi(e_j)\rangle \dv
=&-\int_{\R^m}\big(\eta''(r)x_j+\eta'(r)\frac{x_j}{r}\big)\langle\Delta\tau(\phi),d\phi(e_j)\rangle\dv\\
&-\int_{\R^m}\eta'(r)\frac{x_ix_j}{r}\langle\Delta\tau(\phi),\nabla_id\phi(e_j)\rangle \dv.
\end{align*}

This leads us to the following equality
\begin{align*}
0=&\int_{\R^m}\eta(r)(\frac{m}{2}-3)|\nabla\tau(\phi)|^2\dv
-2\int_{\R^m}\eta'(r)\frac{x_j}{r}\langle\nabla_j\tau(\phi),\tau(\phi)\rangle\dv \\
&+\frac{1}{2}\int_{\R^m}\eta'(r)r|\nabla\tau(\phi)|^2\dv 
-\int_{\R^m}\eta'(r)\frac{x_ix_j}{r}\langle\nabla_i\tau(\phi),\nabla_j\tau(\phi)\rangle\dv \\
&+(m-3)\int_{\R^m}\eta'(r)\frac{x_j}{r}\langle\Delta\tau(\phi),d\phi(e_j)\rangle\dv 
-\int_{\R^m}\eta''(r)x_j\langle\Delta\tau(\phi),d\phi(e_j)\rangle\dv \\
&-2\int_{\R^m}\eta'(r)\frac{x_ix_j}{r}\langle\Delta\tau(\phi),\nabla_id\phi(e_j)\rangle \dv.
\end{align*}

We manipulate the last three terms on the right hand side as follows
\begin{align*}
\int_{\R^m}\eta'(r)\frac{x_j}{r}\langle\Delta\tau(\phi),d\phi(e_j)\rangle\dv
=&\int_{\R^m}\big(\eta''(r)\frac{x_ix_j}{r^2}-\eta'(r)\frac{x_ix_j}{r^3}\big)\langle\nabla_i\tau(\phi),d\phi(e_j)\rangle\dv \\
&+\int_{\R^m}\eta'(r)\frac{1}{r}\langle\nabla_j\tau(\phi),d\phi(e_j)\rangle\dv \\
&+\int_{\R^m}\eta'(r)\frac{x_j}{r}\langle\nabla_i\tau(\phi),\nabla_id\phi(e_j)\rangle\dv
\end{align*}
and
\begin{align*}
\int_{\R^m}\eta''(r)x_j\langle\Delta\tau(\phi),d\phi(e_j)\rangle\dv=
&\int_{\R^m}\eta^{(3)}(r)\frac{x_jx_i}{r}\langle\nabla_i\tau(\phi),d\phi(e_j)\rangle\dv \\
&+\int_{\R^m}\eta''(r)\langle\nabla_j\tau(\phi),d\phi(e_j)\rangle\dv \\
&+\int_{\R^m}\eta''(r)x_j\langle\nabla_i\tau(\phi),\nabla_id\phi(e_j)\rangle\dv.
\end{align*}

Moreover, we find
\begin{align*}
\int_{\R^m}\eta'(r)\frac{x_ix_j}{r}\langle\Delta\tau(\phi),\nabla_id\phi(e_j)\rangle \dv=
&\int_{\R^m}\big(\eta''(r)\frac{x_kx_ix_j}{r^2}-\eta'(r)\frac{x_ix_jx_k}{r^3}\big)\langle\nabla_k\tau(\phi),\nabla_i d\phi(e_j)\rangle\dv \\
&+2\int_{\R^m}\eta'(r)\frac{x_j}{r}\langle\nabla_i\tau(\phi),\nabla_id\phi(e_j)\rangle \dv \\
&+\int_{\R^m}\eta'(r)\frac{x_ix_j}{r}\langle\nabla_k\tau(\phi),\nabla_k\nabla_id\phi(e_j)\rangle \dv.
\end{align*}

Combining the previous equations we find
\begin{align*}
(\frac{m}{2}-3)\int_{\R^m}\eta(r)|\nabla\tau(\phi)|^2\dv
=&2\int_{\R^m}\eta'(r)\frac{x_j}{r}\langle\nabla_j\tau(\phi),\tau(\phi)\rangle\dv
-\frac{1}{2}\int_{\R^m}\eta'(r)r|\nabla\tau(\phi)|^2\dv \\
&+\int_{\R^m}\eta'(r)\frac{x_ix_j}{r}\langle\nabla_i\tau(\phi),\nabla_j\tau(\phi)\rangle\dv\\
&+(3-m)\int_{\R^m}\big(\eta''(r)\frac{x_ix_j}{r^2}-\eta'(r)\frac{x_ix_j}{r^3}\big)\langle\nabla_i\tau(\phi),d\phi(e_j)\rangle\dv \\
&+(3-m)\int_{\R^m}\eta'(r)\frac{1}{r}\langle\nabla_j\tau(\phi),d\phi(e_j)\rangle\dv \\
&+(3-m)\int_{\R^m}\eta'(r)\frac{x_j}{r}\langle\nabla_i\tau(\phi),\nabla_id\phi(e_j)\rangle\dv \\
&+\int_{\R^m}\eta^{(3)}(r)\frac{x_jx_i}{r}\langle\nabla_i\tau(\phi),d\phi(e_j)\rangle\dv \\
&+\int_{\R^m}\eta''(r)\langle\nabla_j\tau(\phi),d\phi(e_j)\rangle\dv \\
&+\int_{\R^m}\eta''(r)x_j\langle\nabla_i\tau(\phi),\nabla_id\phi(e_j)\rangle\dv\\
&+2\int_{\R^m}\big(\eta''(r)\frac{x_kx_ix_j}{r^2}-\eta'(r)\frac{x_ix_jx_k}{r^3}\big)\langle\nabla_k\tau(\phi),\nabla_i d\phi(e_j)\rangle\dv \\
&+4\int_{\R^m}\eta'(r)\frac{x_j}{r}\langle\nabla_i\tau(\phi),\nabla_id\phi(e_j)\rangle \dv \\
&+2\int_{\R^m}\eta'(r)\frac{x_ix_j}{r}\langle\nabla_k\tau(\phi),\nabla_k\nabla_id\phi(e_j)\rangle \dv.
\end{align*}

In order to estimate all terms on the right hand side we first of all note that
\begin{align*}
\big|\int_{\R^m}\eta'(r)r|\nabla\tau(\phi)|^2\dv\big| &
\leq C\int_{B_{2R}\setminus B_R}|\nabla\tau(\phi)|^2\dv,\\
\big|\int_{\R^m}\eta'(r)\frac{x_ix_j}{r}\langle\nabla_i\tau(\phi),\nabla_j\tau(\phi)\rangle\dv\big|&
\leq C\int_{B_{2R}\setminus B_R}|\nabla\tau(\phi)|^2\dv,\\
\big|\int_{\R^m}\eta'(r)\frac{x_ix_j}{r}\langle\nabla_k\tau(\phi),\nabla_k\nabla_id\phi(e_j)\rangle\dv\big| &
\leq C\int_{B_{2R}\setminus B_R}|\nabla^2 d\phi|^2\dv.
\end{align*}

Moreover, we apply the estimates
\begin{align*}
\int_{\R^m}\eta'(r)\frac{x_j}{r}\langle\nabla_j\tau(\phi),\tau(\phi)\rangle\dv &
\leq\frac{C}{R}\|\tau(\phi)\|_{L^2(\R^m)}\|\nabla\tau(\phi)\|_{L^2(\R^m)},\\
\int_{\R^m}\eta''(r)\frac{x_ix_j}{r^2}\langle\nabla_i\tau(\phi),d\phi(e_j)\rangle\dv &
\leq\frac{C}{R^2}\|d\phi\|_{L^2(\R^m)}\|\nabla\tau(\phi)\|_{L^2(\R^m)}, \\
\int_{\R^m}\eta'(r)\frac{x_ix_j}{r^3}\langle\nabla_i\tau(\phi),d\phi(e_j)\rangle\dv &
\leq\frac{C}{R^2}\|d\phi\|_{L^2(\R^m)}\|\nabla\tau(\phi)\|_{L^2(\R^m)}, \\
\int_{\R^m}\eta'(r)\frac{1}{r}\langle\nabla_j\tau(\phi),d\phi(e_j)\rangle\dv & 
\leq\frac{C}{R^2}\|d\phi\|_{L^2(\R^m)}\|\nabla\tau(\phi)\|_{L^2(\R^m)}, \\
\int_{\R^m}\eta'(r)\frac{x_j}{r}\langle\nabla_i\tau(\phi),\nabla_id\phi(e_j)\rangle\dv &
\leq\frac{C}{R}\|\nabla d\phi\|_{L^2(\R^m)}\|\nabla\tau(\phi)\|_{L^2(\R^m)}, \\
\int_{\R^m}\eta^{(3)}(r)\frac{x_jx_i}{r}\langle\nabla_i\tau(\phi),d\phi(e_j)\rangle\dv &
\leq\frac{C}{R^2}\|d\phi\|_{L^2(\R^m)}\|\nabla\tau(\phi)\|_{L^2(\R^m)}, \\
\int_{\R^m}\eta''(r)\langle\nabla_j\tau(\phi),d\phi(e_j)\rangle\dv &
\leq\frac{C}{R^2}\|d\phi\|_{L^2(\R^m)}\|\nabla\tau(\phi)\|_{L^2(\R^m)}, \\
\int_{\R^m}\eta''(r)x_j\langle\nabla_i\tau(\phi),\nabla_id\phi(e_j)\rangle\dv &
\leq\frac{C}{R}\|\nabla d\phi\|_{L^2(\R^m)}\|\nabla\tau(\phi)\|_{L^2(\R^m)}, \\
\int_{\R^m}\eta''(r)\frac{x_kx_ix_j}{r^2}\langle\nabla_k\tau(\phi),\nabla_i d\phi(e_j)\rangle\dv &
\leq\frac{C}{R}\|\nabla d\phi\|_{L^2(\R^m)}\|\nabla\tau(\phi)\|_{L^2(\R^m)}, \\
\int_{\R^m}\eta'(r)\frac{x_ix_jx_k}{r^3}\langle\nabla_k\tau(\phi),\nabla_i d\phi(e_j)\rangle\dv &
\leq\frac{C}{R}\|\nabla d\phi\|_{L^2(\R^m)}\|\nabla\tau(\phi)\|_{L^2(\R^m)}, \\
\int_{\R^m}\eta'(r)\frac{x_j}{r}\langle\nabla_i\tau(\phi),\nabla_id\phi(e_j)\rangle \dv &
\leq\frac{C}{R}\|\nabla d\phi\|_{L^2(\R^m)}\|\nabla\tau(\phi)\|_{L^2(\R^m)}.
\end{align*}

We may conclude that for \(m\neq 6\) the following
inequality holds
\begin{align*}
\int_{B_R}|\nabla\tau(\phi)|^2\dv\leq&\frac{C}{R}\|\nabla d\phi\|_{L^2(\R^m)}\|\nabla\tau(\phi)\|_{L^2(\R^m)}
+\frac{C}{R^2}\|d\phi\|_{L^2(\R^m)}\|\nabla\tau(\phi)\|_{L^2(\R^m)} \\
&+C\int_{B_{2R}\setminus B_R}|\nabla^2 d\phi|^2\dv.
\end{align*}

Taking the limit \(R\to\infty\) 
and using the finiteness assumption \eqref{thm-finiteness-triharmonic}
we conclude that \(\nabla\tau(\phi)=0\).
In order to draw the conclusion that \(\phi\) must be harmonic we calculate
\begin{align*}
0=\int_{\R^m}\eta^2\langle\underbrace{\nabla_i\tau(\phi)}_{=0},d\phi(e_i)\rangle\dv=-\int_{\R^m}\eta^2|\tau(\phi)|^2\dv-2\int_{\R^m}\langle\tau(\phi),d\phi(e_i)\rangle\eta\nabla_i\eta\dv.
\end{align*}
Now, we obtain
\begin{align*}
\frac{1}{2}\int_{\R^m}\eta^2|\tau(\phi)|^2\dv\leq\frac{C}{R^2}\int_{\R^m}|d\phi|^2\dv\to 0 \textrm{ as } R\to\infty
\end{align*}
due to the finiteness assumption \eqref{thm-finiteness-triharmonic}.
The statement that \(\phi\) must be trivial if \(m\neq 2\) follows
from a classical result of Sealey \cite{MR654088}.
The proof is complete.
\end{proof}

\begin{Bem}
Note that Theorem \ref{liouville-triharmonic} and its proof would also hold
if we would consider a weak solution of the triharmonic map equation \eqref{triharmonic-map}.
In this case a weak solution would correspond to \(\phi\in W^{3,2}(\R^m,N)\) that solves
\eqref{triharmonic-map} in a distributional sense.
\end{Bem}

\par\medskip
\textbf{Acknowledgements:}
The author gratefully acknowledges the support of the Austrian Science Fund (FWF) 
through the project P30749-N35 ``Geometric variational problems from string theory''.

\bibliographystyle{plain}
\bibliography{mybib}
\end{document}